\documentclass[11pt]{article}

\usepackage{graphicx}
\usepackage{float}

\usepackage{latexsym}
\usepackage{amssymb}
\usepackage{amsmath}
\usepackage{amsfonts}
\usepackage[ansinew]{inputenc}
\usepackage{tikz}

\def\scat#1#2%
 {\hbox{\vrule\vbox to#2pc{\hrule\hbox to#1pc{\hfil} 
  \vfil\hrule}\vrule}}

\def\findimlogo
 {\scat{.25}{.5}}

\def\findim
 {\hfil\findimlogo\break\par}

\newtheorem{theorem}{Theorem}[section]
\newtheorem{lemma}{Lemma}[section]

\newtheorem{proposition}{Proposition}[section]

\newtheorem{conjecture}{Conjecture}[section]

\newtheorem{myexample}{Example}[section]
\newenvironment{example}
{\begin{myexample}\rm}{\end{myexample}}
\newtheorem{myremark}{Remark}[section]

\newcommand{\pf}{{\bf Proof.~}} 

\newcommand{\lemref}[1]{Lemma~\ref{#1}}

\newcommand{\proposref}[1]{Proposition~\ref{#1}}

\newcommand{\conjref}[1]{Conjecture~\ref{#1}}

\newcommand\Sp
 {\mathop{\rm Sp}\nolimits}

\newcommand\lin
 {\mathop{\rm lin}\nolimits}

\newcommand\aper
 {\mathop{\rm int\,}\nolimits}

\newcommand\conv
 {\mathop{\rm conv\,}\nolimits}

\newcommand\rk
 {\mathop{\rm rank\,}\nolimits}

\newcommand\Hom
 {\mathop{\rm Hom\,}\nolimits}

\newcommand\ord
 {\mathop{\rm ord\,}\nolimits}

\newcommand\card
 {\mathop{\#\,}\nolimits}

\newcommand\re
 {\mathop{\rm Re\,}\nolimits}

\newcommand\im
 {\mathop{\rm Im\,}\nolimits}

\newcommand\Ch
 {\mathop{\rm Ch\,}\nolimits}

\newcommand\Log
 {\mathop{\rm Log\,}\nolimits}

\newcommand\grad
 {\mathop{\rm grad}\nolimits}

\newcommand\Arg
 {\mathop{\rm Arg}\nolimits}
 
\newcommand\Ex
 {\mathop{\rm Exp}\nolimits}
 
\newcommand\Vertx
 {\mathop{\rm Vert}\nolimits}

\begin{document}

\title{The Ronkin number of an exponential sum.}

\author{James Silipo\thanks{The author was partially supported by the grant n. A1UNICAL053 - POR Calabria 2000-2006 Misura 3.7 Azione B.}}

\date{October 29, 2010}

\maketitle

\begin{abstract}
We give an intrinsic estimate of the number of connected components of the complementary set to the amoeba of an exponential sum with real spectrum improving the result of Forsberg, Passare and Tsikh in the polynomial case and that of Ronkin in the exponential one.

\end{abstract}

\section{Introduction}

Let $\Omega\subseteq \mathbb R^n$ be a non empty open convex subset and let $T_\Omega=\Omega+i\mathbb R^n\subseteq \mathbb C^n$ be the (vertical) tube domain on the base $\Omega$. Let also $\re:T_\Omega\rightarrow\Omega$ be the projection of $T_\Omega$ onto its base. The \textit{amoeba} ${\mathcal F}_Y$ of a closed analytic subset $Y$ of $T_\Omega$ is the topological closure in $\Omega$ of $\re Y$, i.e.
\begin{eqnarray}\label{ea}
{\mathcal F}_Y=\overline{\re Y}\,.
\end{eqnarray}

This notion of amoeba was originally proposed by Favorov~\cite{Fav} for zero sets of holomorphic almost periodic functions defined on tube domains.

Recall that a holomorphic function $f$ defined on a tube domain $T_\Omega\subseteq\mathbb C^n$ is said to be \textit{almost periodic} if, for every $D\Subset \Omega$, $f$ is the uniform limit on $T_D=D+i\mathbb R^n$ of a sequence of exponential sums, i.e. $\mathbb C$-linear combinations of exponentials~$e^{\langle z,\lambda\rangle}$, with~$z\in\mathbb C^n$,~$\lambda\in\mathbb R^n$ and~$\langle z,\lambda\rangle=z_1\lambda_1+\ldots+z_n\lambda_n$.
For any non empty open and convex subset $\Omega\subseteq\mathbb R^n$, the class of holomorphic almost periodic functions on $T_\Omega$ will be noted by $HAP(T_\Omega)$, the subclass of exponential sums will be noted by $\Ex (T_\Omega)$.

Any $f\in HAP(T_\Omega)$ has a well defined \textit{Bohr transform} $a(f,\,\cdot\,):\mathbb R^n\rightarrow \mathbb C$ given, for every $\lambda \in\mathbb R^n$, by
\begin{eqnarray}
a(f,\lambda)
=
\lim_{s\to+\infty}
\frac{1}{(2s)^n}
\int_{{\vert y_\ell\vert<s},{\ell=1,\ldots,n}}
e^{-\langle x+iy,\lambda\rangle}
f(x+iy){\rm d}y\,.
\end{eqnarray}
The Bohr transform of $f$ is zero on the whole $\mathbb R^n$ with the exception of a countable set denoted by $\Sp f$ and called the \textit{spectrum} of $f$,
\begin{eqnarray}
\Sp f=\{\lambda\in\mathbb R^n\mid a(f,\lambda)\neq0\}\,.
\end{eqnarray}
The convex hull of $\Sp f$ is noted by $\Gamma_f$, whereas $\Xi_f$ and $\Xi_f^+$ will respectively denote the additive subgroup and the additive submonoid of~$\mathbb R^n$ generated by $\Sp f$. Sometimes we will also consider the linear span $\lin\Gamma_f$ of $\Gamma_f$ in $\mathbb R^n$. When $f$ is an exponential sum, the convex subset $\Gamma_f$ is a polytope which will be referred to as the \textit{Newton polytope} of $f$.

Given a tube domain $T_\Omega$, an important, though very special, subclass of holomorphic almost periodic functions is the class $HSE(T_\Omega)$ provided by those holomorphic almost periodic functions which can be represented as a composition
\begin{eqnarray}
F\circ\psi:T_\Omega\stackrel{\psi}{\longrightarrow}W\stackrel{F}{\longrightarrow}\mathbb C\,,
\end{eqnarray}
where $F$ is a holomorphic function on a multicircular and logarithmically convex domain $W\subseteq \mathbb C^r$ and $\psi$ is a holomorphic mapping of the form
\begin{eqnarray}\label{psi}
\psi(z)=\big(e^{\langle z,\,\omega_1\rangle},\ldots,e^{\langle z,\,\omega_r\rangle}\big),
\end{eqnarray}
with $\omega_1,\ldots,\omega_r\subset\mathbb R^n$ linearly independent over $\mathbb Z$. If $W=(\mathbb C^*)^r$ and $F$ is a Laurent polynomial, one immediately realizes that $\Ex(T_\Omega)\subset HSE(T_\Omega)$.

If $f\in HAP(T_\Omega)$, the corresponding hypersurface amoeba is simply denoted by~${\mathcal F}_f$ and has an interesting concavity property. In fact let $E$ be a connected component of the complementary set $\Omega\setminus{\mathcal F}_f$. The function $f^{-1}$ is holomorphic on the tube $T_E$ and it cannot admit any holomorphic continuation to a strictly larger domain, so by standard facts about holomorphicity domains it follows that $E$ is convex.

In his pioneering article~\cite{Ron} Ronkin proved that for any $f\in HSE(T_\Omega)$ the number of connected components of $\Omega\setminus {\mathcal F}_f$ is locally finite. Favorov~\cite{Fav}, who was the first to use the term ``amoeba'' in the almost periodic context, generalized Ronkin's result to the larger class of holomorphic almost periodic functions with spectrum in a free group.

For $f\in\Ex(T_{\mathbb R^n})$, Ronkin~\cite{Ron} showed that the number $\rho(f)$ of connected components of $\mathbb R^n\setminus {\mathcal F}_f$ is finite and satisfies the estimate 
\begin{eqnarray}\label{ronkinest}
\rho(f)\leqslant\upsilon(f,\gamma)
=
2^{-r}\varkappa_r\big(\sqrt{r}+2r\max_{k\in\gamma(\Sp f)}\Vert k\Vert\big)^r\,,
\end{eqnarray}
where $r=\rk\Xi_f$, $\varkappa_r$ is the Lebesgue measure of the $r$-dimensional unit ball, $\gamma$ is an isomorphism of $\Xi_f$ on $\mathbb Z^r$ and $\Vert\cdot\Vert$ is the usual max norm on $\mathbb R^r$. 
Moreover, if $J_f$ denotes the \textit{Jessen function} $J_f$ of $f$, i.e.
\begin{eqnarray}\label{jessen}
J_f(x)=
\lim_{t\to+\infty}
\frac{1}{(2t)^n}
\int_{[-t,t]^n}
\ln\vert f(x+iy)\vert\, dy
\,,
\end{eqnarray}
Favorov~\cite{Fav} proved that the gradient mapping of $J_f$ is constant on each connected component of $\mathbb R^n\setminus {\mathcal F}_f$ and maps injectively the set of such components into the group $\Xi_f$. 
Stronger results can be obtained in the case of exponential sums with integer spectrum, i.e. $\Sp f\subset\mathbb Z^n$. In fact in this case the amoeba theory is equivalent to the \textit{polynomial amoeba theory}. In order to introduce this polynomial approach, consider the proper mapping
\begin{eqnarray}
\Log:(\mathbb C^*)^r\longrightarrow\mathbb R^r
\end{eqnarray}
that maps a point $\zeta=(\zeta_1,\ldots,\zeta_r)\in(\mathbb C^*)^r$ to the point $(\ln\vert \zeta_1\vert,\ldots,\ln\vert \zeta_r\vert)$. 
Following Passare and Tsikh~\cite{PT}, given a closed analytic subset $X\subset(\mathbb C^*)^r$, the amoeba of $X$ is the image ${\mathcal A}_X$ of $X$ under the mapping $\Log$, i.e.
\begin{eqnarray}\label{pa}
{\mathcal A}_X=\Log (X)\,.
\end{eqnarray}

If $X\subset (\mathbb C^*)^r\subset\mathbb C^r$ is a closed analytic subset, then ${\mathcal A}_X\neq{\mathcal F}_X$, so in order to distinguish these clearly different notions, an amoeba of the form~(\ref{pa}) will be referred to as a \textit{polynomial amoeba} whereas one of the form~(\ref{ea}) will be called an \textit{exponential amoeba}.

The mapping~(\ref{psi}) with $r=n$ and $\omega_1,\ldots,\omega_n$ equal to the canonical basis of $\mathbb R^n$ yields the relation
\begin{eqnarray}
{\mathcal A}_X={\mathcal F}_{\psi^{-1}(X)}\,,
\end{eqnarray}
i.e. the polynomial theory of amoebae is a special issue of the exponential one. If $P\in \mathbb C[\zeta_1^{\pm 1},\ldots,\zeta_r^{\pm 1}]$ is a non constant Laurent polynomial, the amoeba of its zero set, often noted by ${\mathcal A}_P$, is a proper closed subset of $\mathbb R^r$ whose complementary set has a finite number of connected components each of which is convex. The number of such components clearly equals the Ronkin number of the exponential sum $f=P\circ\psi$ and for sake of simplicity this number will be denoted $\rho(P)$ instead of $\rho(P\circ\psi)$. Forsberg, Passare and Tsikh~\cite{FPT} have shown that, for any Laurent polynomial $P\in \mathbb C[\zeta_1^{\pm 1},\ldots,\zeta_r^{\pm 1}]$, the following estimate holds true
\begin{eqnarray}\label{fopats}
\card\Vertx\Gamma_P\leqslant\rho(P)\leqslant\card(\Gamma_P\cap\mathbb Z^r)\,,
\end{eqnarray}
where $\Vertx\Gamma_P$ is the set of vertices of the Newton polytope $\Gamma_P$ of $P$ and $\Gamma_P\cap\mathbb Z^r$ is the set of lattice points belonging to $\Gamma_P$, (cf.Gelfand et al.~\cite{GKZ}, Forsberg et al.~\cite{FPT}). An injective mapping from the set of connected components of $\mathbb R^n\setminus{\mathcal A}_P$ into $\Gamma_P\cap\mathbb Z^r$ is provided by the gradient of the so-called \textit{Ronkin function} of $P$, i.e. the convex function defined for every $x\in\mathbb R^r$ as
\begin{eqnarray}\label{ronkf}
N_P(x)
=
\int_{\Log^{-1}(x)}
\ln\vert P(\zeta)\vert\,d\eta_r(\zeta)\,,
\end{eqnarray}
where $\eta_r(\zeta)$ is the unique translation invariant probability Haar measure on the real $r$-dimensional torus~$\Log^{-1}(x)$. This measure can be computed through the differential form
\begin{eqnarray}
\frac{1}{(2\pi i)^r}\frac{d\zeta_1\,\wedge\ldots\,\wedge d\zeta_r}{\zeta_1\cdot\ldots\cdot\zeta_r}
\end{eqnarray}
or equivalently through the form
\begin{eqnarray}
\frac{1}{(2\pi)^r}\,d\Arg\zeta_1\wedge\,\ldots,\,\wedge \,d\Arg\zeta_r\,,
\end{eqnarray}
on~$[0,2\pi]^r$. In Number Theory, the Ronkin function of a Laurent polynomial $P$ is known as the \textit{Mahler measure} of $P$, it is well defined even on the amoeba~${\mathcal A}_P$, where the integrand is manifestly singular. On the amoeba complement the function $N_P$ is piecewise linear and its gradient defines the orders of the components of $\mathbb R^r\setminus{\mathcal A}_P$, (cf. Forsberg et al.~\cite{FPT}). In particular, for a complement component~$E\subset\mathbb R^r\setminus{\mathcal A}_P$, the \textit{order}~$\ord(E)$ is defined as the value of~$\grad N_P(x)$ for any $x\in E$. Different components have different orders and the order of a complement component always belongs to $\Gamma_P\cap\mathbb Z^r$. 
For a thorough exposition of these facts the reader is referred to Forsberg et al.~\cite{FPT} and Rullg{\aa}rd~\cite{Ru1}.

It should be noted that for an exponential sum with integer spectrum the estimate~(\ref{fopats}) is sharper than~(\ref{ronkinest}). Our main result improves both these estimates. It can be summarized as follows.

\smallskip
\noindent
Let~$f\in \Ex(\mathbb C^n)$ with $\rk\Xi_f=r$ and let~$\gamma$ be a group isomorphism of~$\Xi_f$ onto~$\mathbb Z^r$. Then
\begin{enumerate}
	\item the set
	\begin{eqnarray}
	\Lambda_f=\gamma^{-1}\big(\conv(\gamma(\Sp f))\cap\gamma(\Xi_f)\big)\,,
	\end{eqnarray}
	where the convex hull~$\conv(\gamma(\Sp f))$ is taken in~$\mathbb R^r=\gamma(\Xi_f)\otimes_\mathbb Z \mathbb R$, is a finite subset of $\Xi_f\cap\Gamma_f$ which contains $\Sp f$ and does not depend on the isomorphism~$\gamma$;
		\item the number $\rho(f)$ satisfies the estimate
	\begin{eqnarray}
\card\Vertx(\Gamma_f)
\leqslant
\rho(f)
\leqslant
\card\Lambda_f\,,
\end{eqnarray}
where $\Vertx(\Gamma_f)$ denotes the set of vertices of $\Gamma_f$;
\item the gradient of $J_f$ injects the set of connected components of $\mathbb R^n\setminus{\mathcal F}_f$ into $\Lambda_f$;
\item  $\card\Lambda_f<\upsilon(f,\gamma)$;
\item if $\Sp f\subset\mathbb Z^n$ and $\Xi_f\subsetneq \mathbb Z^n$, then $\card\Lambda_f<\card(\mathbb Z^n\cap \Gamma_f)$. 
\end{enumerate}

In the sequel of this article, for any $f\in\Ex(\mathbb C^n)$, the number $\rho(f)$ of connected components of $\mathbb R^n\setminus{\mathcal F}_f$ will be referred to as \textit{the Ronkin number} of $f$. Analogously, for any Laurent polynomial $P$ the number $\rho(P)$ of connected components of $\mathbb R^n\setminus{\mathcal A}_P$ will be referred to as \textit{the Ronkin number} of $P$.

\section{Counting components}

Let $\Ch\mathbb R^n=\Hom_\mathbb Z(\mathbb R^n,\mathbb S^1)$ be the (multiplicative, abelian) group of~$\mathbb S^1$-characters of~$\mathbb R^n$. 
If~$f\in\Ex(\mathbb C^n)$ and~$\chi\in\Ch\mathbb R^n$, the exponential sum
\begin{eqnarray}
f_\chi(z)
=
\sum_{\lambda\in\Sp f}
a(f,\lambda)\,\chi(\lambda)\,e^{\langle z,\lambda\rangle}
\end{eqnarray}
is called the \textit{perturbation} of~$f$ by~$\chi$. This yields an action
\begin{eqnarray}
\Ch\mathbb R^n\times\Ex(\mathbb C^n)\longrightarrow \Ex(\mathbb C^n)
\end{eqnarray}
such that~$(fg)_\chi=f_\chi g_\chi$, for any exponential sums~$f,g$ and any character~$\chi$. 

Given~$f\in\Ex(\mathbb C^n)$, the orbit of $f$ for this action is merely the set of perturbations of~$f$ by the characters in the group~$\Ch\Xi_f=\Hom_\mathbb Z(\Xi_f,\mathbb S^1)$. 

Further details on such perturbations can be found in Silipo~\cite{sil}, or in Fabiano et al.~\cite{FGS} for a generalization. 

If~$f\in\Ex(\mathbb C^n)$ and~$\chi\in\Ch\mathbb R^n$, then~$V(f_\chi)$ denotes the zero set of $f_\chi$. The following result, (a proof of which is available in Silipo~\cite{sil} Proposition~3.2 and Corollary~3.9, or in Fabiano et al.~\cite{FGS} Theorem~4.3 and Corollary~4.4) will be useful in the sequel.
\begin{theorem}\label{cavallo}
Let $f$ be an exponential sum on~$\mathbb C^n$, then
\begin{eqnarray}\label{cav}
{\mathcal F}_f
=
\mathbb R^n
\cap
\bigcup_{\chi\in\Ch\Xi_f}
V(f_\chi)
=
\bigcup_{\chi\in\Ch\Xi_f}
\re V(f_\chi)
\,,
\end{eqnarray}
in particular~${\mathcal F}_f={\mathcal F}_{f_\chi}$, for every~$\chi\in\Ch\Xi_f$. 
\end{theorem}

\medskip
\noindent
Another easy though useful fact about the shape of an exponential amoeba is given by the following proposition.

\begin{proposition}\label{fulldim}
Let~$f$ be an exponential sum defined on~$\mathbb C^n$, then
\begin{eqnarray}
{\mathcal F}_f
=
{\mathcal F}_f+(\lin\Gamma_f)^\perp\,,
\end{eqnarray}
where $(\lin\Gamma_f)^\perp$ is the orthogonal complement of the linear subspace spanned by $\Gamma_f$ in  $\mathbb R^n$ endowed with the standard scalar product.
\end{proposition}
\pf
Let~$\chi\in\Ch\Xi_f$ and consider the zero set~$V(f_\chi)$ of~$f_\chi$ in~$\mathbb C^n$. If~$\pi_f$ is the linear projection of~$\mathbb C^n$ onto the complex linear subspace~$(\lin\Gamma_f+i\lin\Gamma_f)$, then, the expression of~$f_\chi$ implies that, 
\begin{eqnarray}
z\in V(f_\chi)\Longleftrightarrow \pi_f(z)+(\lin\Gamma_f+i\lin\Gamma_f)^{\perp_\mathbb C}\subseteq V(f_\chi)\,,
\end{eqnarray}
where $\perp_\mathbb C$ stands for orthogonality in~$\mathbb C^n$ with its standard hermitian product. Taking real parts yields
\begin{eqnarray}
z\in V(f_\chi)\Longrightarrow \re \pi_f(z)+(\lin\Gamma_f)^\perp\subseteq \re V(f_\chi)\,,
\end{eqnarray}
where the equality~$\re (\lin\Gamma_f+i\lin\Gamma_f)^{\perp_\mathbb C}=(\lin\Gamma_f)^\perp$ follows by the obvious inclusion~$\Gamma_f\subset\mathbb R^n$. As a consequence
\begin{eqnarray}
\bigcup_{z\in V(f_\chi)}
\re \pi_f(z)+ (\lin\Gamma_f)^\perp
\subseteq 
\re V(f_\chi)\,,
\end{eqnarray}
and taking the union on~$\Ch\Xi_f$ implies that
\begin{eqnarray*}
\bigcup_{\chi\in\Ch\Xi_f}
\bigcup_{z\in V(f_\chi)}
\re \pi_f(z)+(\lin\Gamma_f)^\perp
&=&
(\lin\Gamma_f)^\perp
+
\bigcup_{\chi\in\Ch\Xi_f}
\bigcup_{z\in V(f_\chi)}
\re \pi_f(z)
\\
&=&
(\lin\Gamma_f)^\perp
+
\bigcup_{\chi\in\Ch\Xi_f}
\re V(f_\chi)
\\
&=&
(\lin\Gamma_f)^\perp
+
{\mathcal F}_f
\,,
\end{eqnarray*}
so~${\mathcal F}_f+(\lin\Gamma_f)^\perp$ has to be included in~${\mathcal F}_f$ thus proving what claimed.\findim

The idea behind our counting technique is to linearly embed an exponential amoeba into a naturally associated polynomial one so as to obtain an intrinsic estimate of the Ronkin number. We start by constructing the associated polynomial amoeba.

Let~$f\in\Ex(\mathbb C^n)$ and let $\gamma:\Xi_f\rightarrow\mathbb Z^{\rk \Xi_f}$ be a group isomorphism. 
In the vector space $\mathbb Z^{\rk \Xi_f}\otimes_\mathbb Z \mathbb R$ consider the convex hull $\conv(\gamma(\Sp f))$ of the image of $\Sp f$ via $\gamma$ and the set 
\begin{eqnarray}
\Lambda_f=\gamma^{-1}(\conv(\gamma(\Sp f))\cap \gamma(\Xi_f))\,.
\end{eqnarray}
Since $\gamma(\Xi_f)=\mathbb Z^{\rk\Xi_f}$, the set $\Lambda_f$ is the subset of $\Xi_f$ consisting of the inverse images via $\gamma$ of the lattice points belonging to the polytope $\conv(\gamma(\Sp f))$. Evidently, among these lattice points there are the elements of $\Sp f$, so that $\Sp f\subseteq\Lambda_f$, but in fact one can say more.

\begin{lemma}\label{conv}
Let~$f\in\Ex(\mathbb C^n)$ with $\rk\Xi_f=r$ and let $\gamma:\Xi_f\rightarrow\mathbb Z^r$ be a group isomorphism. Then $\Lambda_f\subset\Gamma_f$.
\end{lemma}
\pf Let $\mu\in\Lambda_f$, then there exists a $t\in [0,1]^{\Sp f}$ summing up to $1$ such that
\begin{eqnarray}
\gamma(\mu)=\sum_{\lambda\in\Sp f} t_\lambda \gamma(\lambda)\in\conv(\gamma(\Sp f))\cap\mathbb Z^r.
\end{eqnarray}
For every $1\leqslant \ell\leqslant r$, let $\omega_\ell$ the inverse image via $\gamma$ of the $\ell$-th element $e_\ell$ of the canonical basis of $\mathbb Z^r$, then, for every $\lambda\in\Sp f$,
\begin{eqnarray}
\lambda=\sum_{\ell=1}^r \lambda_\ell\,\omega_\ell\,,\qquad
\gamma(\lambda)=\sum_{\ell=1}^r \lambda_\ell\,e_\ell
\end{eqnarray}
and
\begin{eqnarray}
\gamma(\mu)=\sum_{\ell=1}^r\bigg(\sum_{\lambda\in\Sp f} t_\lambda \lambda_\ell\bigg) e_\ell\in \mathbb Z^r\,,
\end{eqnarray}
so that
\begin{eqnarray}
\mu
&=&
\gamma^{-1}
\left(
\sum_{\ell=1}^r\bigg(\sum_{\lambda\in\Sp f} t_\lambda \lambda_\ell\bigg) e_\ell
\right)
\\
&=&
\sum_{\ell=1}^r\bigg(\sum_{\lambda\in\Sp f} t_\lambda\lambda_\ell\bigg)\omega_\ell
\\
&=&
\sum_{\lambda\in\Sp f} t_\lambda\bigg(\sum_{\ell=1}^r\lambda_\ell\,\omega_\ell\bigg)
\\
&=&
\sum_{\lambda\in\Sp f} t_\lambda \lambda\in\Gamma_f\,.
\end{eqnarray}
By the arbitrary choice of $\mu$ in $\Lambda_f$ we obtain $\Lambda_f\subset\Gamma_f$.\findim

The following lemma implies that $\Lambda_f$ does not depend on the isomorphism used in its definition.

\begin{lemma}\label{comb}
Let $f$ be an exponential sum and let~$\gamma_1$ and $\gamma_2$ be two different isomorphisms of the group~$\Xi_f$ on $\mathbb Z^{\rk\Xi_f}$. Then the polytopes $\conv(\gamma_1(\Sp f))$ and $\conv(\gamma_2(\Sp f))$ are combinatorially isomorphic.
\end{lemma}
\pf Let $r=\rk \Xi_f$. The authomorphism of~$\mathbb Z^r$ given by~$\gamma_1\circ\gamma_2^{-1}$ has an obvious continuation to an $\mathbb R$-linear automorphism of~$\mathbb R^r$. The two lattice polytopes $\conv(\gamma_1(\Sp f))$ and $\conv(\gamma_2(\Sp f))$ correspond to each other in this $\mathbb R$-linear authomorphism, hence they are combinatorially isomorphic. \findim
\par
\medskip

Since a combinatorial isomorphism of lattice polytopes preserves lattice points in the polytopes, the set $\Lambda_f$ proves to depend only on $f$. As shown in the sequel of the article, the interest in the set $\Lambda_f$ is due to the role of ``order set'' it will play for $\mathbb R^n\setminus{\mathcal F}_f$. 

Now, if $r=\rk\Xi_f$, fix an isomorphism $\gamma:\Xi_f\to\mathbb Z^r$ and for every $1\leqslant j\leqslant r$, let $\omega_j$ be the inverse image via $\gamma$ of the $j$-th element of the canonical basis of $\mathbb Z^r$. This basis and the isomorphism $\gamma$ will be referred to as \textit{associated} to each other. Then consider the Laurent polynomial
\begin{eqnarray}\label{polinomio}
P(\zeta)=\sum_{\lambda\in\Sp f}a(f,\lambda)\zeta^{\gamma(\lambda)}
=
\sum_{k\in\gamma(\Sp f)}
a(f,\gamma^{-1}(k))\,
\zeta_1^{k_1}\ldots\,\zeta_r^{k_r}
\,,
\end{eqnarray}
and observe that for any character $\chi\in\Ch\mathbb R^n$ we have
\begin{eqnarray}\label{import}
f_\chi(x)
=
P\big(e^{\langle x,\,\omega_1\rangle+i\Arg\chi(\omega_1)},\ldots,e^{\langle x,\,\omega_r\rangle+i\Arg\chi(\omega_r)}\big)\,,
\end{eqnarray}
for every $x\in\mathbb R^n$. Consider also the linear mapping~$L:\mathbb R^n\rightarrow\mathbb R^r$ given by
\begin{equation}\label{lo}
L(x)=(\langle x,\omega_1\rangle,\ldots,\langle x,\omega_r\rangle)\,,
\end{equation}
and notice that the kernel of~$L$ equals the orthogonal complement (with respect to the standard scalar product of $\mathbb R^n$) of the linear subspace $\lin\Gamma_f$ spanned by~$\Gamma_f$. It follows that
\begin{eqnarray}
\dim(L(\mathbb R^n))
=
n-\dim\ker L
=
\dim\lin\Gamma_f
\,,
\end{eqnarray}
in particular, the mapping~$L$ is injective if and only if the linear subspace~$\lin\Gamma_f$ is full-dimensional.
\par
\medskip
The following lemma gives some additional information about this construction.
\begin{lemma}\label{ok}
Let $f\in\Ex(\mathbb C^n)$ and $\gamma:\Xi_f\to\mathbb Z^r$ be an isomorphism. Let also $P$ be the corresponding Laurent polynomial,~$\,\Sigma_P$ the normal fan of~$\Gamma_P$,~$h_{\Gamma_P}$ the support function of the polytope $\Gamma_P$ and, for every face ~$\Delta$ of~$\Gamma_P$, let~$K_{\Delta,\Gamma_P}$ be the corresponding dual cone. Then
\begin{enumerate}
\item $L(\mathbb R^n)\setminus\{0\}$ does not intersect the cones of~$\,\Sigma_P$ corresponding to the positive dimensional faces of~$\Gamma_P$,
\item \begin{eqnarray}
h_{\Gamma_P}(L(x))
=
h_{\Gamma_f}(x)\,,
\end{eqnarray}
for every~$x\in\mathbb R^n$. In particular, for every~$\lambda\in\Vertx(\Gamma_f)$, 
\begin{eqnarray}
L(K_{\lambda,\Gamma_f})\subseteq K_{\gamma(\lambda),\Gamma_P}\,.
\end{eqnarray}
\item \begin{eqnarray}
L(\mathbb R^n)
\subset
\{0\}\cup
\left(
\bigcup_{\lambda\in\Vertx(\Gamma_f)}
\aper K_{\gamma(\lambda),\Gamma_P}
\right)\,.
\end{eqnarray}
\end{enumerate}
\end{lemma}
\pf 1. Let~$\{\omega_1,\ldots,\omega_r\}$ be the basis of~$\,\Xi_f$ associated to $\gamma$ and suppose, by contradiction, there is an~$x\in\mathbb R^n\setminus\ker L$ such that~$L(x)$ belongs to a cone of~$\Sigma_P$ corresponding to a positive dimensional face~$\Delta$ of~$\Gamma_P$. Since~$\Delta$ is positive dimensional, it admits two distinct lattice points~$u$ and~$v$, so the vectors~$t=u-v$ and~$L(x)$ are orthogonal to each other. This means that
\begin{eqnarray}
0=
\langle t,L(x)\rangle_{\mathbb R^r}
=
\sum_{\ell=1}^r
t_\ell\,
\langle x,\omega_\ell\rangle_{\mathbb R^n} 
=
\sum_{\ell=1}^r
\langle x,t_\ell\,\omega_\ell\rangle_{\mathbb R^n} 
=
\langle
x,
\sum_{\ell=1}^r
t_\ell\,\omega_\ell
\rangle_{\mathbb R^n}\,.
\end{eqnarray}
Since~$\omega_1,\ldots,\omega_r$ are $\mathbb Z$-linear independent and~$t\neq 0$, the preceding equalities make~$x$ to belong to~$(\lin \Gamma_f)^\perp$, which equals the kernel of~$L$. The contradiction implies the statement.

2. Let~$x\in\mathbb R^n$ be fixed, then by definition of support function
\begin{eqnarray}
h_{\Gamma_P}(L(x))
=
\sup_{v\in\Gamma_P}\langle L(x),v\rangle_{\mathbb R^r}
=
\sup_{v\in\Vertx(\Gamma_P)}\langle L(x),v\rangle_{\mathbb R^r}\,.
\end{eqnarray}
However~$\Vertx(\Gamma_P)\subseteq\gamma(\Sp f)$, so
\begin{eqnarray}
h_{\Gamma_P}(L(x))
=
\sup_{\lambda\in\Sp f}
\langle L(x),\gamma(\lambda)\rangle_{\mathbb R^r}\,.
\end{eqnarray}
Any~$\lambda\in\Sp f$ has an expression of the form
\begin{eqnarray}
\lambda
=
\sum_{\ell=1}^r
k_{\lambda,\ell}\,\omega_\ell\,,
\end{eqnarray}
for a uniquely determined sequence $k_{\lambda,1},\ldots,k_{\lambda,r}$ of integers, then
\begin{eqnarray}
h_{\Gamma_P}(L(x))
=
\sup_{\lambda\in\Sp f}
\sum_{\ell=1}^r
\langle x,\omega_\ell\rangle_{\mathbb R^n}
k_{\lambda,\ell}
=
\sup_{\lambda\in\Sp f}
\langle x,
\sum_{\ell=1}^r
k_{\lambda,\ell}\,\omega_\ell
\rangle_{\mathbb R^n}
=
\sup_{\lambda\in\Sp f}
\langle x,\lambda\rangle_{\mathbb R^n}\,,
\end{eqnarray}
i.e.~$h_{\Gamma_P}(L(x))=h_{\Gamma_f}(x)$. 

In particular, if~$\lambda\in\Gamma_f$ and~$x\in K_{\lambda,\Gamma_f}$, then
\begin{eqnarray}
\langle L(x),\gamma(\lambda)\rangle_{\mathbb R^r}=\langle x,\lambda\rangle_{\mathbb R^n}=h_{\Gamma_f}(x)=h_{\Gamma_P}(L(x))\,,
\end{eqnarray}
i.e.~$L(K_{\lambda,\Gamma_f})\subseteq K_{\gamma(\lambda),\Gamma_P}$.

3. The whole $\mathbb R^n$ is equal to the union of the closures of the cones which are dual to the vertices of~$\Gamma_f$ and each of these cones is mapped by~$L$ in the interior of the dual cone associated to the corresponding vertex in~$\Gamma_P$.\findim
\par
\medskip

\begin{theorem}\label{io}
Let $f\in\Ex(\mathbb C^n)$ and let $\gamma:\Xi_f\rightarrow\mathbb Z^r$ be an isomorphism with the associated basis $\{\omega_1,\ldots,\omega_r\}$ of~$\,\Xi_f$. If $P$ is the corresponding Laurent polynomial and $L$ the corresponding linear mapping, then
\begin{enumerate}
\item the following two equalities hold true \begin{eqnarray}
L({\mathcal F}_f)=L(\mathbb R^n)\cap {\mathcal A}_P\qquad {\it and}\qquad
L(\mathbb R^n\setminus{\mathcal F_f})=L(\mathbb R^n)\setminus {\mathcal A}_P\,,
\end{eqnarray}
\item the Ronkin number~$\rho(f)$ equals the number of connected components of~$L(\mathbb R^n)\setminus{\mathcal A}_P$,
\item the following estimate holds true 
\begin{eqnarray}\label{grezza}
\card\Vertx(\Gamma_f)
\leqslant
\rho(f)
\leqslant
\card\Lambda_f\,.
\end{eqnarray}
\end{enumerate}
\end{theorem}
\pf
1. Let us start with the first equality. If~$x\in{\mathcal F}_f$, then, by equality~(\ref{cav}), there exists a character~$\chi\in\Ch\Xi_f$ such that~$f_\chi(x)=0$. By~(\ref{import}) it follows that~$L(x)\in{\mathcal A}_P$. Conversely, if~$x\in\mathbb R^n$ is such that~$L(x)$ belongs to~${\mathcal A}_P$, then there is a zero~$\zeta$ of~$P$ for which~$\Log(\zeta)=L(x)$, i.e.~$\ln\vert\zeta_\ell\vert=\langle x,\omega_\ell\rangle$, for~$1\leqslant \ell\leqslant r$. Let~$\chi\in\Ch\Xi_f$ be the uniquely defined character such that~$\chi(\omega_\ell)=e^{i\im \zeta_\ell}$, for~$1\leqslant \ell\leqslant r$, then 
\begin{eqnarray}
0=P(\zeta)=f_\chi(x)\,,
\end{eqnarray}
or equivalently,~$x\in{\mathcal F}_f$.

As for the second equality, if $x$ does not belong to~${\mathcal F}_f$, then~$f_\chi(x)\neq0$, for any character. Suppose, by contradiction, that $L(x)$ belongs to~${\mathcal A}_P$. Then~$L(x)=\Log(\zeta)$, for some zero $\zeta$ of $P$. Let~$\vartheta\in[0,2\pi)^r$ be the $r$-tuple of principal arguments of $\zeta$ and let also~$\chi\in\Ch\Xi_f$ be the corresponding character. By virtue of~(\ref{import}),~$f_\chi(x)=0$ and this contradicts the choice of~$x$, thus~$L(\mathbb R^n\setminus{\mathcal F}_f)\subseteq L(\mathbb R^n)\setminus{\mathcal A}_P$.
Finally, let~$x\in\mathbb R^n$ such that~$L(x)\notin{\mathcal A}_P$, then for every~$\vartheta\in[0,2\pi)^r$, 
\begin{eqnarray}
P\big(e^{\langle x,\,\omega_1\rangle+i\vartheta_1},\ldots,e^{\langle x,\,\omega_r\rangle+i\vartheta_r}\big)\neq 0\,,
\end{eqnarray}
equivalently, (again by~(\ref{import})),~$f_\chi(x)\neq0$ for any~$\chi\in\Ch\Xi_f$, whence~$x\notin{\mathcal F}_f$.

2. The connected components of~$L(\mathbb R^n)\setminus{\mathcal A}_P$ are convex. Let~$m$ be the number of such components. Since $L$ is continuous, by~\lemref{io}, we know that~$\rho(f)\geqslant m$. Suppose, by contradiction, $\rho(f)>m$. Then we may find two points~$x_1$ and~$x_2$ belonging to distinct components of~$\mathbb R^n\setminus{\mathcal F}_f$ which are mapped by~$L$ to a same component~$Y$ of~$L(\mathbb R^n)\setminus{\mathcal A}_P$. By virtue of~\proposref{fulldim}, the projections~$x_1^\prime$ and~$x_2^\prime$ on~$\lin\Gamma_f$ of~$x_1$ and~$x_2$ respectively are still distinct and belonging to~${\mathcal F}_f$. Now,~$\ker L=(\lin\Gamma_f)^\perp$, so the restriction of~$L$ to~${\mathcal F}_f\cap\lin\Gamma_f$ realizes a linear homeomorphism onto~$L(\mathbb R^n)\setminus{\mathcal A}_P$. If~$\alpha\subset Y$ is a line segment joining~$L(x_1)$ and~$L(x_2)$, then~$L^{-1}(\alpha)\cap\lin\Gamma_f$ is a line segment joining~$x_1^\prime$ and~$x_2^\prime$, thus~\proposref{fulldim} makes~$x_1$ and~$x_2$ to belong to the same component of~$\mathbb R^n\setminus{\mathcal F}_f$. The contradiction implies the corollary.

3. The results of \cite{FPT} imply that the number of components of $\mathbb R^r\setminus{\mathcal A}_P$ cannot exceed $\card(\Gamma_P\cap\mathbb Z^r)$. Since the linear subspace~$L(\mathbb R^n)$ may intersect all the components of~$\mathbb R^r\setminus{\mathcal A}_P$, by \lemref{comb} we get the desired estimate from above. The lower bound is easily found by the usual geometric series trick.  \findim

It should perhaps be mentioned that the estimate~(\ref{grezza}) completely neglects the values of the coefficients of $f$ except to those which correspond to the vertices of $\Gamma_f$, on which the only requirement is to be non-zero. A better estimate would involve in a more substantial way the values of the Fourier coefficients of $f$. Nevertheless, we notice that when $\Sp f\subseteq \mathbb Z^n$, the present estimate improves the well known result of Forsberg, Passare and Tsikh \cite{FPT} since the upper bound $\card\Lambda_f$ does not exceed $\card(\Gamma_f\cap\mathbb Z^n)$. As an example, consider the exponential sum $f(z)=2+e^{2z_1}+e^{2z_2}+e^{4z_1+4z_2}\in\Ex(\mathbb C^2)$. Of course $\Xi_f=(2\mathbb Z)^2$, so if $\gamma:\Xi_f\rightarrow\mathbb Z^2$ is the isomorphism mapping $(2,0)$ to $(1,0)$ and $(0,2)$ to $(0,1)$, it follows that $\card(\Lambda_f)=5$, whereas $\card(\Gamma_f\cap\mathbb Z^2)=11$.

\section{Some remarks}
\subsection{Maximal sparseness}

A Laurent polynomial is called \textit{maximally sparse} if all the points in the support of summation of the polynomial are vertices of its Newton polytope. 

Likewise, an exponential sum is called \textit{maximally sparse} if all the points in its spectrum are vertices of its Newton polytope.

A group isomorphism $\gamma:\Xi_f\rightarrow\mathbb Z^{\rk\Xi_f}$ cannot generally admit an $\mathbb R$-linear continuation to $\mathbb R^n$, nevertheless, though such an isomorphism cannot preserve the convex structure of $\Gamma_f$, something still survives. The following proposition is a key ingredient in order to relate the notion of maximal sparseness with the preceding construction.

\begin{proposition}\label{vertici}
Let~$f\in\Ex(\mathbb C^n)$ and let $\{\omega_1,\ldots,\omega_r\}$ be the basis of~$\,\Xi_f$ associated to some isomorphism~$\gamma:\Xi_f\rightarrow\mathbb Z^r$ together with the corresponding Laurent polynomial $P$. Then
\begin{eqnarray}
\gamma(\Vertx(\Gamma_f))\subseteq\Vertx(\Gamma_P)\,.
\end{eqnarray}
In particular, if $f$ is maximally sparse, then the corresponding Laurent polynomial~$P$ is maximally sparse too.
\end{proposition}
\pf Suppose, by contradiction, there is a~$\lambda_*\in\Vertx(\Gamma_f)$ such that $\gamma(\lambda_*)$ is not a vertex of~$\Gamma_P$. Since~$\Gamma_P=\conv(\gamma(\Sp f))$, it follows that
\begin{eqnarray}
\gamma(\lambda_*)
=
\sum_{\lambda\in\Sp f\setminus\{\lambda_*\}}
t_\lambda\gamma(\lambda)\,,
\end{eqnarray}
for some family~$\{t_\lambda\}_{\lambda\in\Sp f\setminus\{\lambda_*\}}$ of non negative real numbers summing up to~$1$. For every~$\lambda\in\Sp f$ there is a unique sequence~$k_{\lambda,1},\ldots,k_{\lambda,r}$ of integers such that
\begin{eqnarray}
\lambda
=
\sum_{j=1}^r
k_{\lambda,j}\,\omega_j\,,
\end{eqnarray}
so
\begin{eqnarray*}
\gamma(\lambda_*)
&=&
\sum_{j=1}^r
k_{\lambda_*,j}\,\gamma(\omega_j)
\end{eqnarray*}
and also
\begin{eqnarray*}
\gamma(\lambda_*)
=
\sum_{\lambda\in\Sp f\setminus\{\lambda_*\}}
t_\lambda
\sum_{j=1}^r
k_{\lambda,j}\,\gamma(\omega_j)
=
\sum_{j=1}^r
\left(
\sum_{\lambda\in\Sp f\setminus\{\lambda_*\}}
t_\lambda
k_{\lambda,j}
\right)
\gamma(\omega_j)\,.
\end{eqnarray*}
As~$\{\gamma(\omega_1),\ldots,\gamma(\omega_r)\}$ is a basis of the vector space~$\mathbb Z^r\otimes\mathbb R$, the two last expressions of~$\gamma(\lambda_*)$ yield the equality
\begin{eqnarray}
k_{\lambda_*,j}
=
\sum_{\lambda\in\Sp f\setminus\{\lambda_*\}}
t_\lambda
k_{\lambda,j}\,,
\end{eqnarray}
for every~$1\leqslant j\leqslant r$. In each of these equalities the first member is an integer and so the second member has to be an integer too, it follows that
\begin{eqnarray*}
\sum_{j=1}^r
k_{\lambda_*,j}\,
\omega_j
&=&
\sum_{j=1}^r
\left(
\sum_{\lambda\in\Sp f\setminus\{\lambda_*\}}
t_\lambda
k_{\lambda,j}
\right)
\omega_j\,,
\end{eqnarray*}
or equivalently
\begin{eqnarray}
\lambda_*
=
\sum_{j=1}^r
\left(
\sum_{\lambda\in\Sp f\setminus\{\lambda_*\}}
t_\lambda
k_{\lambda,j}
\right)
\omega_j
=
\sum_{\lambda\in\Sp f\setminus\{\lambda_*\}}
t_\lambda
\sum_{j=1}^r
k_{\lambda,j}\,
\omega_j
=
\sum_{\lambda\in\Sp f\setminus\{\lambda_*\}}
t_\lambda
\lambda\,,
\end{eqnarray}
i.e. $\lambda_*$ is not a vertex of~$\Gamma_f$. The contradiction implies the first statement.

To prove the second statement, notice that for any $f$, $\Vertx(\Gamma_P)\subseteq\gamma(\Sp f)$. If $f$ is maximally sparse, then~$\Sp f=\Vertx(\Gamma_f)$ and
\begin{eqnarray}
\gamma(\Sp f)
=
\gamma(\Vertx(\Gamma_f))
=
\Vertx(\Gamma_P)\,.
\end{eqnarray}
As~$\gamma(\Sp f)$ is precisely the support of summation of~$P$ the conclusion is that all the elements in the support of summation of~$P$ are vertices of~$\Gamma_P$, i.e. $P$ is maximally sparse.
\findim

Observe that the second statement in \proposref{vertici} may not be reversed. In fact the exponential sum $f=1+e^z+e^{\sqrt 2 z}$ is not maximally sparse, whereas~$P=1+\zeta_1+\zeta_2$ is such.

A polynomial (resp. exponential) amoeba is said to be \textit{solid} if its complementary set has the minimal number of connected components.   

Equivalently a polynomial (resp. exponential) amoeba is solid if the number of connected components of its complementary set equals the number of vertices in the Newton polytope of a Laurent polynomial (resp. exponential sum) which defines the amoeba. 

Passare and Rullg{\aa}rd~\cite{PR2} suggested the following conjecture.

\begin{conjecture}[\cite{PR2}]\label{pass-rull}
A maximally sparse Laurent polynomial has a solid amoeba.
\end{conjecture}

The conjecture is true if the Newton polytope of the given maximally sparse Laurent polynomial is reduced to a line segment, however it is not known if the conjecture is true in the general case. Nisse~\cite{N} has recently proposed a solution in the affirmative, but his proof seems to need some clarification. We notice here that such a solution would imply an exponential counterpart.

\begin{proposition}
If \conjref{pass-rull} is true, a maximally sparse exponential sum has a solid amoeba.
\end{proposition}
\pf Let $f$ be a maximally sparse exponential sum on~$\mathbb C^n$. With the same notation of \proposref{vertici}, the polynomial~$P$ is maximally sparse too so by \conjref{pass-rull} the amoeba~${\mathcal A}_P$ is solid.
The estimate~(\ref{grezza}) for maximally sparse~$f$ and~$P$ becomes 
\begin{eqnarray}\label{prima-stima}
\card\Vertx (\Gamma_f)\leqslant \rho(f)\leqslant\card\Vertx(\Gamma_P)\,, 
\end{eqnarray}
but the proof of~\proposref{vertici} shows that the bounds in the above estimate coincide so ${\mathcal F}_f$ is solid.\findim

\subsection{Order theory}

Let $f\in\Ex(\mathbb C^n)$ and let $E$ be a component of~$\;\mathbb R^n\setminus{\mathcal F}_f$. According to Favorov~\cite{Fav}, the {\textit{order}} of~$E$ is the value $\ord(E)$ taken by the gradient of the Jessen function~(\ref{jessen}) of~$f$ at any point of the domain~$E$, i.e.
\begin{eqnarray}
\ord(E)
=
\grad J_f(x)\,,
\end{eqnarray}
for any $x\in E$. Observe that the order of a complement component does not depend merely on the amoeba. In fact, for every $\lambda\in\mathbb R^n$, the exponential sum $g(z)=e^{\langle z,\lambda\rangle} f(z)$ has the same amoeba as $f$ but, for every complement component $E$ and any $x\in E$, one has 
\begin{eqnarray}
\grad J_g(x)=\lambda+\grad J_f(x)\,.
\end{eqnarray}

In the almost periodic literature the order of a component $E\subset \mathbb R^n\setminus{\mathcal F}_f$ is known as the \textit{mean motion} of $f$ in the domain $E$ and, as shown in Ronkin~\cite{Ron}, it is well defined. Here we propose a different approach to this order theory. 
Consider the Pontriagin group~$\Ch\Xi_f=\Hom_{\mathbb Z}(\Xi_f,\mathbb S^1)$. If~$\rk \Xi_f=r$, then~$\Ch\Xi_f$ is an~$r$-dimensional real compact torus. Define the \textit{Ronkin function} of $f$ as the function~${\mathcal N}_f:\mathbb R^n\longrightarrow\mathbb R$ given, for any~$x\in\mathbb R^n$, by
\begin{eqnarray}\label{ronkf2}
{\mathcal N}_f(x)
=
\int_{\Ch\Xi_f}
\ln\vert f_\chi(x)\vert\,d\vartheta(\chi)\,,
\end{eqnarray}
where~$\vartheta(\chi)$ is the translation invariant probability Haar measure on $\Ch\Xi_f$.
The measure~$\vartheta$ can be computed by choosing a basis~$\{\omega_1,\ldots,\omega_r\}$ of~$\,\Xi_f$ and integrating the differential form
\begin{eqnarray}
\frac{1}{(2\pi)^r}
\,
d\Arg\chi(\omega_1)\,\wedge\,\ldots\,\wedge\,d\Arg\chi(\omega_r)\,.
\end{eqnarray}

In order to realize that this definition is a good one, it is enough to consider the Laurent polynomial~$P$ corresponding to a chosen basis of $\Xi_f$ and remark that, for any~$x\in\mathbb R^n$,
\begin{eqnarray*}
{\mathcal N}_f(x)
&=&
\int_{\Ch\Xi_f}
\ln
\vert 
P\big(e^{\langle x,\omega_1\rangle+i\Arg\chi(\omega_1)},\ldots,e^{\langle x,\omega_r\rangle+i\Arg\chi(\omega_r)}\big)
\vert
\,
d\vartheta(\chi)
\\
&=&
\int_{\Log^{-1}(L(x))}
\ln
\vert 
P(\zeta)
\vert
\,
d\eta_r(\zeta)
\,,
\end{eqnarray*}
i.e.
\begin{eqnarray}\label{utile}
{\mathcal N}_f(x)
=
N_P(L(x))\,.
\end{eqnarray}
Thus the Ronkin function of~$f$ is the restriction of the Ronkin function of~$P$ to the subspace~$L(\mathbb R^n)$. As such~${\mathcal N}_f$ proves to be well defined and convex.

The function~${\mathcal N}_f$ was already studied by Ronkin in~\cite{Ron} where he proved the following theorem.
\begin{theorem}[\cite{Ron}, Theorem 6] Let $f$ be an exponential sum, then
\begin{eqnarray}
J_f={\mathcal N}_f\,.
\end{eqnarray}
\end{theorem}

Using the function~${\mathcal N}_f$ and the equality~(\ref{utile}) it is quite easy to prove the following result.

\begin{lemma}\label{ronklem}
Let~$f\in\Ex(\mathbb C^n)$, then the Ronkin function~${\mathcal N}_f$ is piecewise linear on~$\mathbb R^n\setminus {\mathcal F}_f$, its gradient mapping~$\grad{\mathcal N}_f$ realizes an injection of the set of components of~$\mathbb R^n\setminus{\mathcal F}_f$ into~$\Lambda_f$.
\end{lemma}
\pf
Let~$\{\omega_1,\ldots,\omega_r\}$ be a basis of~$\,\Xi_f$ and let~$\gamma$ be the associated isomorphism of~$\Xi_f$ on~$\mathbb Z^r$. If~$P$ is the corresponding Laurent polynomial, then for every~$1\leqslant j\leqslant n$ and any~$x\in \mathbb R^n\setminus{\mathcal F}_f$,
\begin{eqnarray}
\frac{\partial {\mathcal N}_f}{\partial x_j}(x)
=
\sum_{\ell=1}^r
\frac{\partial N_{P}}{\partial \zeta_\ell}\big(L(x)\big)
\frac{\partial (L)_\ell}{\partial x_j}(x)
=
\sum_{\ell=1}^r
\frac{\partial N_{P}}{\partial \zeta_\ell}\big(L(x)\big)
\omega_{\ell,j}\,,
\end{eqnarray}
i.e.
\begin{eqnarray}\label{ordine}
\grad_x{\mathcal N}_f(x)
=
\sum_{\ell=1}^r
\frac{\partial N_P}{\partial \zeta_\ell}(L(x))\,
\omega_\ell
=
\gamma^{-1}(\grad_\zeta N_P(L(x))\,.
\end{eqnarray}
Now,~$N_P$ is piecewise linear on~$\mathbb R^r\setminus{\mathcal A}_{P}$, its gradient mapping is constant on each component of~$\mathbb R^r\setminus{\mathcal A}_P$ and it maps injectively the set of these components on a finite subset of~$\Gamma_P\cap\mathbb Z^r$. Consequently,~${\mathcal N}_f$ is piecewise linear on~$\mathbb R^n\setminus{\mathcal F}_f$, its gradient mapping is constant on each component of~$\mathbb R^n\setminus{\mathcal F}_f$ and, as the~$\omega_1,\ldots,\omega_r$ are~$\mathbb Z$-linearly independent, it maps injectively the set of such components on a finite subset of~$\mathbb R^n$. In order to identify this subset observe that equation~(\ref{ordine}) implies
\begin{eqnarray}
\grad_x{\mathcal N}_f(x)
\in
\gamma^{-1}(\grad_\zeta N_P(L(\mathbb R^n)\setminus{\mathcal A}_P))\subseteq\gamma^{-1}(\Gamma_P\cap\mathbb Z^r)=\Lambda_f\,,
\end{eqnarray}
for any~$x\in \mathbb R^n\setminus{\mathcal F}_f$.\findim

With respect to Favorov's result, \lemref{ronklem} actually adds that, for any exponential sum $f$, the gradient of~$J_f={\mathcal N}_f$ maps the amoeba complement~${\mathbb R}^n\setminus{\mathcal F}_f$ in the Newton polytope~$\Gamma_f$, thus completing the analogy with the polynomial case. 

It should also be noticed that, for a given $f\in\Ex(\mathbb C^n)$, if $r={\rk \Xi_f}$ and $\gamma:\Xi_f\rightarrow\mathbb Z^r$ is a fixed isomorphism then, up to a multiplication by a suitable exponential monomial, the corresponding Laurent polynomial $P$ is an ordinary polynomial and the estimate $\rho(f)\leqslant \upsilon(f,\gamma)$ provided by Ronkin can be found by bounding the number of lattice points (with positive coordinates) belonging to the ball about the origin of $\mathbb R^r$ with radius equal to
\begin{eqnarray}
\max_{x\in\mathbb R^r\setminus{\mathcal A}_P}
\grad N_P(x)\,.
\end{eqnarray}
As the Newton polytope $\Gamma_P$ is properly contained in that ball, \lemref{ronklem} implies that 
\begin{eqnarray}
\card\Lambda_f< \upsilon(f,\gamma)\,.
\end{eqnarray}
It is also worth noting that, unlike $\upsilon(f,\gamma)$, the bound $\card(\Lambda_f)$ is intrinsic since it does not depend on the isomorphism used to compute it. 

\section{Some examples}

The following examples show the preceding constructions at work. For the sake of simplicity we just consider the case $r=2$ and~$1\leqslant n\leqslant r$. The figures placed at the end show, for each example, the Newton polytope $\Gamma_P$, the points in the subset $\gamma(\Lambda_f)\subset\Gamma_P$ and the amoeba ${\mathcal A}_P$ cut by the linear subspace $L(\mathbb R^n)$.

\begin{example}\label{es1}
{\rm Let~$f\in\Ex(\mathbb C)$ be given by~$f(z)=1+3\,e^{\sqrt{2} z}+e^{\sqrt{5} z}$. The group~$\Xi_f$ has rank~$r=2$, the numbers~$\sqrt2$ and~$\sqrt5$ generate it freely and if $\gamma$ denotes the associated isomorphism of $\Xi_f$ on $\mathbb Z^2$ one gets 
\begin{description}
\item[-] $P(\zeta_1,\zeta_2)=1+3\,\zeta_1+\zeta_2\,$,
\item[-] $23<\upsilon(f,\gamma)=2^{-2}\pi\big(\sqrt{2}+4)^2< 24\,$,
\item[-] $\card\Lambda_f= 3\,$.
\end{description}
}
\end{example}

\begin{figure}[H]
\centering
\framebox{
\begin{tikzpicture}[scale=2]
\draw[help lines] (-0.5cm,2.5cm) grid (1.5cm,4.5cm);
\draw (0,3) -- (1,3) -- (0,4) -- (0,3);
\filldraw [black] 
(0,3) circle (0.5pt)
(1,3) circle (0.5pt) 
(0,4) circle (0.5pt) ;
\draw(0.5,2.5) node {\scriptsize $\Gamma_P$};
\end{tikzpicture}
\qquad
\includegraphics[width=0.4\textwidth]{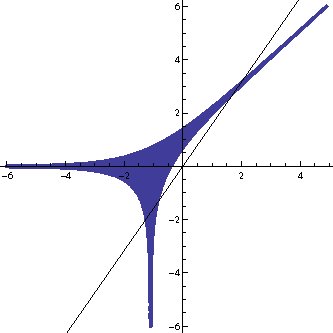}
}
\caption{$\rho(f)=\card\Lambda_f=3$.}
	\label{fig1}
\end{figure}
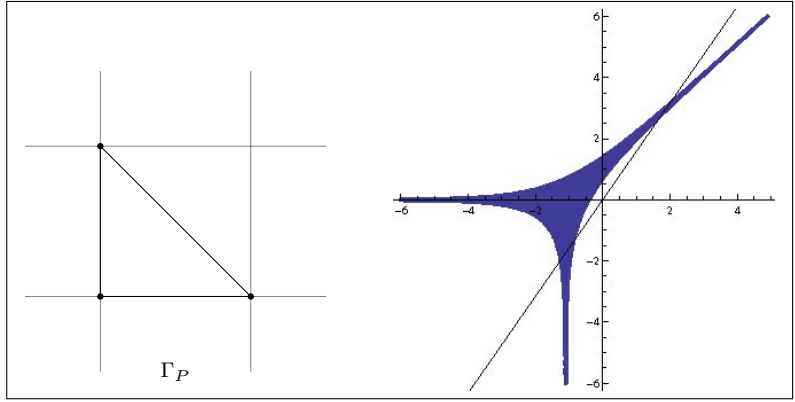

\begin{example}\label{es2}
{\rm Let~$f\in\Ex(\mathbb C)$ be given by~$f(z)=1+e^{3\sqrt2 z}+e^{3\pi z}-6e^{(\sqrt2+\pi)z}$. The group~$\Xi_f$ has rank~$r=2$, the numbers~$3\sqrt2$ and~$\sqrt2+\pi$ generate it freely and if $\gamma$ denotes the associated isomorphism of $\Xi_f$ on $\mathbb Z^2$ one obtains 
\begin{description}
\item[-] $P(\zeta_1,\zeta_2)=1+\zeta_1+\zeta_1^{-1}\zeta_2^3-6\zeta_2\,$,
\item[-] $141<\upsilon(f,\gamma)=2^{-2}\pi\big(\sqrt{2}+12)^2< 142\,$,
\item[-] $\card\Lambda_f= 4\,$.
\end{description}
}
\end{example}

\begin{figure}[H]
\centering
\framebox{
\begin{tikzpicture}[scale=1.2]
\draw[help lines] (-1.5cm,1.75cm) grid (1.5cm,5.5cm);
\draw (0,2) -- (1,2) -- (-1,5) -- (0,2);
\filldraw [black] 
(0,2) circle (1pt)
(1,2) circle (1pt)
(-1,5) circle (1pt)
(0,3) circle (1pt);
\draw(0,1.5) node {\scriptsize $\Gamma_P$};
\end{tikzpicture}
\qquad
\includegraphics[width=0.4\textwidth]{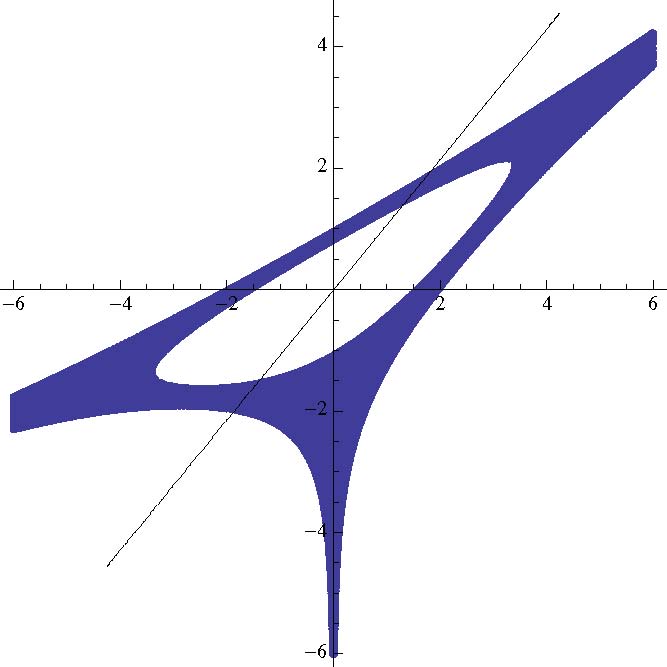}

}
\caption{$3=\rho(f)\leqslant\card\Lambda_f=4$.}
	\label{fig2}
\end{figure}

\begin{example}\label{es3}
{\rm Let~$f\in\Ex(\mathbb C)$ be given by $f(z)=1 + 9e^{\sqrt2 z} + 9 e^{2\sqrt2 z} + e^{3\sqrt2 z} - 9 e^{\sqrt7 z} + 9 e^{2\sqrt7 z} + e^{3\sqrt7 z} - 5e^{(\sqrt2+2\sqrt7) z} - 5 e^{(2\sqrt2+\sqrt7) z} - 35e^{(\sqrt2+\sqrt7) z}$. The group~$\Xi_f$ has rank~$r=2$, the numbers~$\sqrt2$ and~$\sqrt7$ generate it freely and if $\gamma$ denotes the associated isomorphism of $\Xi_f$ onto $\mathbb Z^2$ one obtains 
\begin{description}
\item[-] $P(\zeta_1,\zeta_2)=1+9\zeta_1+9\zeta_1^2+\zeta_1^3-9\zeta_2+9\zeta_2^2+\zeta_2^3-5\zeta_1\zeta_2^2-5\zeta_1^2\zeta_2-35\zeta_1\zeta_2\,$,
\item[-] $141<\upsilon(f,\gamma)=2^{-2}\pi\big(\sqrt{2}+12)^2< 142\,$,
\item[-] $\card\Lambda_f=10\,$.
\end{description}
}
\end{example}

\begin{figure}[H]
\centering
\framebox{
\begin{tikzpicture}
\draw[help lines] (-0.5cm,1.5cm) grid (3.5cm,5.5cm);
\draw (0,2) -- (3,2) -- (0,5) -- (0,2);
\filldraw [black] 
(0,2) circle (1pt)
(1,2) circle (1pt)
(2,2) circle (1pt)
(3,2) circle (1pt)
(0,3) circle (1pt)
(0,4) circle (1pt)
(0,5) circle (1pt)
(1,3) circle (1pt)
(1,4) circle (1pt)
(2,3) circle (1pt)
(0,3) circle (1pt);
\draw(1.5,1.5) node {\scriptsize $\Gamma_P$};
\end{tikzpicture}
\qquad
\includegraphics[width=0.4\textwidth]{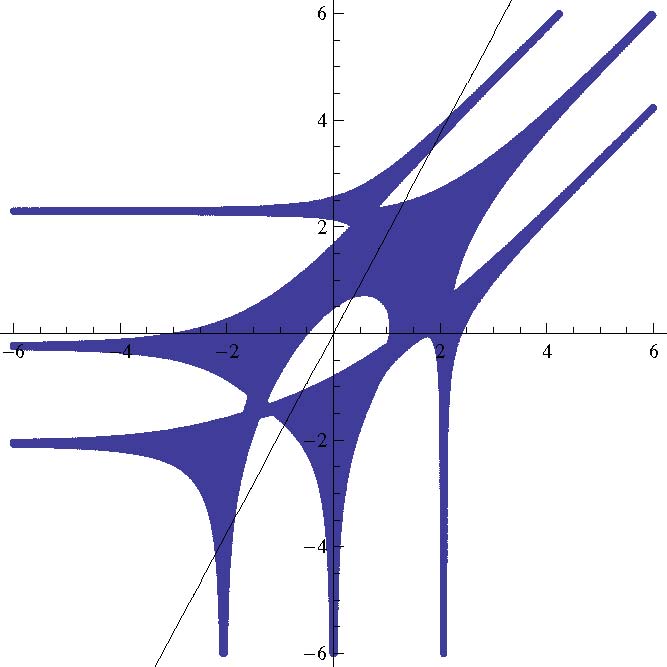}
}
\caption{$5=\rho(f)\leqslant\card\Lambda_f=10$.}
	\label{fig3}
\end{figure}

\begin{example}\label{es4}
{\rm Let~$f\in\Ex(\mathbb C)$ be given by $f(z)=2 + 9 e^{\sqrt2 z} + 9 e^{\sqrt7 z} + e^{2\sqrt2 z} + e^{2\sqrt7 z} + 180 e^{(\sqrt2+\sqrt7)z} + 9e^{(2\sqrt2+\sqrt7)z} + 9e^{(\sqrt2+2\sqrt7)z} + e^{2(\sqrt2+\sqrt7)z}$. The group~$\Xi_f$ has rank~$r=2$, the numbers~$\sqrt2$ and~$\sqrt7$ generate it freely and if $\gamma$ denotes the associated isomorphism of $\Xi_f$ on $\mathbb Z^2$ one gets 

\begin{description}
\item[-] $P(\zeta_1,\zeta_2)=2+9\zeta_1+9\zeta_2+\zeta_1^2+\zeta_2^2+180\zeta_1\zeta_2+9\zeta_1^2\zeta_2+9\zeta_1\zeta_2^2+\zeta_1^2\zeta_2^2\,$,
\item[-] $69<\upsilon(f,\gamma)=2^{-2}\pi\big(\sqrt{2}+8)^2< 70\,$,
\item[-] $\card\Lambda_f=9\,$.
\end{description}
}
\end{example}

\begin{figure}[H]
\centering
\framebox{
\begin{tikzpicture}
\draw[help lines] (-0.5cm,1.5cm) grid (2.5cm,4.5cm);
\draw (0,2) -- (2,2) -- (2,4) -- (0,4) -- (0,2);
\filldraw [black] 
(0,2) circle (1pt)
(1,2) circle (1pt)
(2,2) circle (1pt)
(2,3) circle (1pt)
(0,3) circle (1pt)
(0,4) circle (1pt)
(1,3) circle (1pt)
(2,4) circle (1pt)
(1,4) circle (1pt);
\draw(1,1) node {\scriptsize $\Gamma_P$};
\end{tikzpicture}
\qquad
\includegraphics[width=0.4\textwidth]{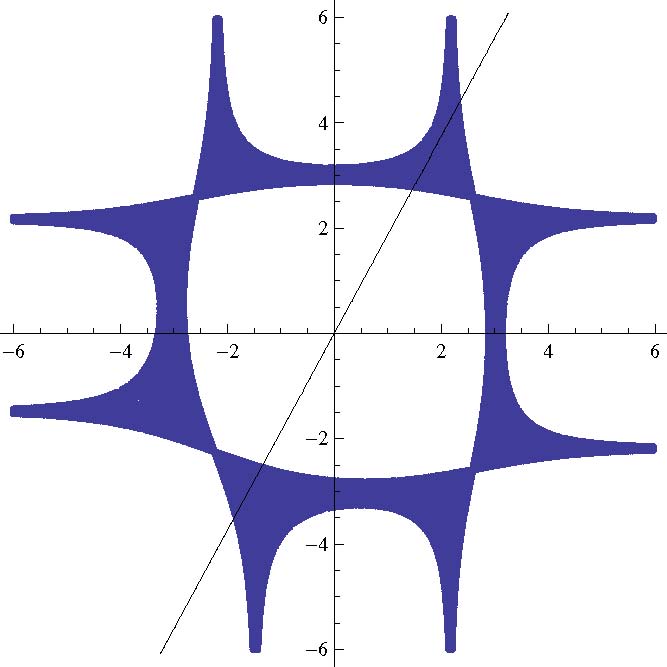}
}
\caption{$3=\rho(f)\leqslant\card\Lambda_f=9$.}
	\label{fig4}
\end{figure}

\begin{example}\label{es5}
{\rm Let~$f\in\Ex(\mathbb C^2)$ be given by $f(z_1,z_2)=2 + 7 e^{2 z_1} + 9 e^{4 z_2} + e^{2 z_1+4 z_2} + 18 e^{z_1+2z_2} $. The group~$\Xi_f$ has rank~$r=2$, the elements~$(2,0)$ and~$(1,2)$ generate it freely and if $\gamma$ denotes the associated isomorphism of $\Xi_f$ on $\mathbb Z^2$ one gets 
\begin{description}
\item[-] $P(\zeta_1,\zeta_2)=2+7\zeta_1+9\zeta_1^{-1}\zeta_2^2+\zeta_2^2+18\zeta_2\,$,
\item[-] $69<\upsilon(f,\gamma)=2^{-2}\pi\big(\sqrt{2}+8)^2< 70 \,$,
\item[-] $\card(\Gamma_f\cap\mathbb Z^2)=15\,$,
\item[-] $\card\Lambda_f=5\,$.
\end{description}
}
\end{example}

\begin{figure}[H]
\centering
\framebox{
\begin{tikzpicture}
\draw[help lines] (-0.5cm,1.5cm) grid (2.5cm,6.5cm);
\draw[thick] (0,2) -- (2,2) -- (2,6) -- (0,6) -- (0,2);
\filldraw [black] 
(0,2) circle (1pt)
(0,3) circle (1pt)
(0,4) circle (1pt)
(0,5) circle (1pt)
(0,6) circle (1pt)
(1,2) circle (1pt)
(1,3) circle (1pt)
(1,4) circle (1pt)
(1,5) circle (1pt)
(1,6) circle (1pt)
(2,2) circle (1pt)
(2,3) circle (1pt)
(2,4) circle (1pt)
(2,5) circle (1pt)
(2,6) circle (1pt);
\draw(1,1) node {\scriptsize $\Gamma_f$};
\end{tikzpicture}
\qquad
\begin{tikzpicture}
\draw[help lines] (-1.5cm,1.5cm) grid (1.5cm,4.5cm);
\draw[thick] (0,2) -- (1,2) -- (0,4) -- (-1,4) -- (0,2);
\filldraw [black] 
(0,2) circle (1pt)
(0,3) circle (1pt)
(0,4) circle (1pt)
(1,2) circle (1pt)
(-1,4) circle (1pt);
\draw(0,1) node {\scriptsize $\Gamma_P$};
\end{tikzpicture}
\qquad
\includegraphics[width=0.4\textwidth]{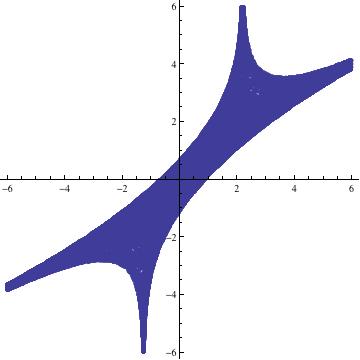}
}
\caption{$4=\card\Vertx\Gamma_f=\rho(f)\leqslant\card\Lambda_f=5$.}
	\label{fig7}
\end{figure}

\end{document}